\documentclass[11pt]{article}
\usepackage{amssymb}
\usepackage{mathrsfs}
\usepackage{latexsym,amsmath,amssymb,amsfonts,epsfig,graphicx,cite,psfrag}
\usepackage{eepic,color,colordvi,amscd}
\usepackage{enumerate}
\usepackage{enumitem}

\newcommand{\qed}{\hfill $\Box $}
\newcommand{\pf}{\noindent {\bf Proof.} }
\newtheorem{theorem}{Theorem}
\newtheorem{lemma}[theorem]{Lemma}

\newtheorem{obser}[theorem]{Observation}
\newtheorem{cla}{Claim}

\newtheorem{conjecture}[theorem]{Conjecture}

\begin{document}

\title{Improved  Bound on Vertex Degree Version of Erd\H{o}s Matching Conjecture}
 
\author{ Mingyang Guo, Hongliang Lu\footnote{Partially  supported by the National Natural
Science Foundation of China under grant No.12271425}  and Yaolin Jiang\footnote{Partially  supported by  the International Science and Technology Cooperation Program of Shaanxi Key Research \& Development Plan under grant 2019KWZ-08}\\
School of Mathematics and Statistics\\
Xi'an Jiaotong University\\
Xi'an, Shaanxi 710049, China
}
\date{}


\date{}

\maketitle

\begin{abstract}
For a $k$-uniform hypergraph $H$, let $\delta_1(H)$ denote the minimum vertex degree of $H$, and $\nu(H)$ denote the
size of the largest  matching in $H$. In this paper, we show that for any $k\geq 3$ and $\beta>0$, there exists an integer $n_0(\beta,k)$ such that for positive integers $n\geq n_0$ and $m\leq (\frac{k}{2(k-1)}-\beta)\frac{n}{k}$, if $H$ is an $n$-vertex $k$-graph with $\delta_1(H)>{{n-1}\choose {k-1}}-{{n-m}\choose {k-1}},$  then  $\nu(H)\geq m$. This improves upon  earlier results of Bollob\'{a}s, Daykin  and Erd\H{o}s (1976) for the range $n> 2k^3(m+1)$ and Huang and Zhao (2017) for the range $n\geq 3k^2 m$.
\end{abstract}

\section{Introduction}
Let $k$ be a positive integer. For a set
$S$, let ${S\choose k}:=\{T\subseteq S: |T|=k\}$. A {\it hypergraph} $H$
consists of a vertex set $V(H)$ and an edge set $E(H)$ whose members
are subsets of $V(H)$. A hypergraph $H$ is {\it $k$-uniform} if $E(H)\subseteq {V(H)\choose k}$, and a $k$-uniform hypergraph is also
called a {\it $k$-graph}. We use $e(H)$ to denote the  number of edges of $H$.

Let $H$ be a $k$-graph and $T\subseteq V(H)$. The {\it degree}
of $T$ in $H$, denoted by $d_H(T)$, is the number of edges of $H$
containing $T$. Let $\ell$ be a nonnegative integer; then $\delta_\ell(H):=\min\{d_H(T): T\in {V(H)\choose \ell}\}$ denotes the
{\it minimum $\ell$-degree} of $H$. Hence,  $\delta_0(H)$ is the number of edges in $H$.  Note that $\delta_1(H)$ is often called the  {\it minimum vertex degree} of $H$, and $\delta_{k-1}(H)$ is also known as the {\it minimum codegree} of $H$.
A {\it matching} in $H$ is a set of pairwise disjoint edges of $H$, and it is  called a \emph{perfect matching} if the union of all edges of the matching is $V(H)$. We use $\nu(H)$ to denote the size of the largest matching in $H$.
Let $K_r^k$ denote a complete $k$-graph with $r$ vertices and edge set ${V(K_r^k)\choose k}$.

Erd\H{o}s and Gallai \cite{EG59} determined the threshold of $\delta_0(G)$ for a 2-graph $G$ to contain a matching of given size, and  Erd\H{o}s~\cite{Er65} conjectured
the following generalization to $k$-graphs for $k\ge 3$: The threshold of $\delta_0(H)$ for a $k$-graph $H$ on $n$ vertices to contain a matching of size $m$ is
$$\max\left\{{km-1\choose k},{n\choose k}-{n-m+1\choose k}\right\}+1.$$
For recent results on this conjecture, we refer the reader to \cite{AFHRRS12,Fr13,FK18,HLS}.

 R\"{o}dl, Ruci\'{n}ski and Szemer\'{e}di
\cite{RRS09} determined the minimum codegree threshold for  the
existence of perfect matchings in $k$-graphs.
It is conjectured in \cite {KO09, HPS09} that the $\ell$-degree threshold for the existence of a perfect matching in a $k$-graph $H$ is
$$\delta_\ell(H)> \left(\max \left\{\frac 1 2, 1-\left(1-\frac 1 k\right)^{k-\ell}\right\}+o(1)\right){n-\ell\choose {k-\ell}}$$
for $k\geq 3$ and $1\leq \ell<k$.
The first term $(1/2+o(1)){n-\ell\choose k-\ell}$ comes from a parity construction: Take disjoint nonempty sets
$A$ and $B$ with $||A|-|B||\le 2$ and $|A|\equiv 1 \pmod 2$, and form a hypergraph $H$ by taking all $k$-subsets $f$ of $A\cup B$ with $|f\cap A|\equiv0 \pmod 2$.
The second term is given by the hypergraph obtained from $K_n^{k}$ (the complete $k$-graph on $n$ vertices) by deleting
all edges from a subgraph isomorphic to  $K_{n-n/k+1}^{k}$.
Treglown and Zhao \cite{TZ12,TZ13}
determined the minimum $\ell$-degree threshold for the appearance of perfect matchings in $k$-graphs, for $k/2\leq \ell\leq k-2$.
For a 3-graph $H$, H\`{a}n, Person and Schacht \cite{HPS09} showed that $\delta_1(H)>(5/9+o(1)){|V(H)|\choose 2}$
is sufficient for the appearance of a perfect matching of  $H$.
 K\"{u}hn, Osthus and Treglown \cite{KOT13} proved a stronger result:
There exists a positive integer $n_0$ such that if $H$ is a 3-graph with $|V(H)|=n\geq n_0$, $m$ is an integer with $1\le m\le n/3$,
and $\delta_1(H)>{{n-1}\choose 2}-{{n-m}\choose 2},$ then $\nu(H)\ge m$. For $k\in \{3,4\}$, Khan \cite{Kh13,Kh16}  showed  that there exists a positive integer $n_0$ such that if $H$ is a $k$-graph with $|V(H)|=n\geq n_0$ and $\delta_1(H)>{{n-1}\choose k-1}-{{n-n/k}\choose k-1},$ then $H$ has a perfect matching, where $n\equiv 0\pmod k$.

%
Bollob\'{a}s, Daykin  and Erd\H{o}s \cite{BDE76}  proved that for integers  $k\geq 2$ and $m\ge 1$, if $H$ is a $k$-graph with $|V(H)|=n> 2k^3(m+1)$  and
$\delta_1(H)>{{n-1}\choose {k-1}}-{{n-m}\choose {k-1}},$
then $\nu(H)\ge m$. The bound is tight. To see it, we define an $n$-vertex $k$-graph $H_{k}(U,W)=(V,E)$, where $V$ is partitioned into $U\cup W$, and  $E=\{e\in {V \choose k} : 1\leq |e\cap W|\leq k-1\}$.
When $|W|=m-1$, it is easy to see that $\delta_1(H_{k}(U,W))={n-1\choose k-1}-{n-m\choose k-1}$ and $\nu(H_{k}(U,W))=m-1$. We denote $H_{k}(U,W)$ with $|U|=n-m+1$ and $|W|=m-1$ by $H_k(n,m)$.

This result was improved by Huang and Zhao \cite{HZ17}, who proved that for $n\geq 3k^2m$, if $\delta_1(H)>{{n-1}\choose {k-1}}-{{n-m}\choose {k-1}}$, then $\nu(H)\ge m$.
Huang and Zhao \cite{HZ17} also proposed the following conjecture.
\begin{conjecture}[Huang and Zhao, \cite{HZ17}]\label{conj1}
	Given positive integers $m,k,n$ such that $m<n/k$, let $H$ be a $k$-graph on $n$ vertices. If $\delta_1(H)>\binom{n-1}{k-1}-\binom{n-m}{k-1}$, then $\nu(H)\geq m$.
\end{conjecture}

Kupavskii \cite{Kupavskii19} proved an $\ell$-degree version result for $1\leq \ell\leq k-1$, which confirmed Conjecture \ref{conj1} for $n\geq 2k^2$ and $k\geq3(m-1)$. When $n\geq 2(k-1)m$, an asymptotic version of this conjecture was proved by K\"uhn, Osthus and Townsend \cite{KOT14}. For fractional matchings, Huang and Zhao \cite{HZ17} proved that when $n\geq(2m-1)(k-1)-m+2$, every $k$-graph $H$ on $n$ vertices with $\delta_1(H)>\binom{n-1}{k-1}-\binom{n-m}{k-1}$ contains a fractional matching of size $m$. Frankl and Kupavskii \cite{FK18} improved the result of K\"uhn, Osthus and Townsend \cite{KOT14} and the fractional version result of Huang and Zhao \cite{HZ17} by extending the range of $n$ to $n\geq \frac{5}{3}(k-1)m-\frac{2}{3}m+1$.

In the paper, we obtain the following result.
\begin{theorem}\label{main-thm}
Let $k\geq 3$ be an integer. For any $\beta>0$, there exists an integer $n_0=n_0(\beta,k)$ such that the following holds. Let $n,m$ be integers such that $n\geq n_0$ and $1\leq m\leq (\frac{k}{2(k-1)}-\beta)\frac{n}{k}$.  Let $H$ be a $k$-graph on $n$ vertices. If $\delta_1(H)>{n-1\choose k-1}-{n-m\choose k-1}$, then $\nu(H)\geq m$.

\end{theorem}

Given two $k$-graphs $H_1, H_2$ and a real number $\varepsilon>0$, we say that $H_2$ \textit{$\varepsilon$-contains} $H_1$ if $V(H_1) = V(H_2)$ and $|E(H_1)\setminus E(H_2)|\leq\varepsilon|V(H_1)|^k$. Let $H$ be an $n$-vertex $k$-graph. We say that $H$ \textit{$\varepsilon$-contains} $H_{k}(n,m)$ if there exists a partition $U,W$ of $V(H)$ with $|W|=m-1$ and $H$ $\varepsilon$-contains $H_{k}(U,W)$. Otherwise, we say that $H$ does not \textit{$\varepsilon$-contain} $H_{k}(n,m)$. Our proof of Theorem~\ref{main-thm} consists of two parts by considering whether or not  $H$ $\varepsilon$-contains $H_{k}(n,m)$, which  is similar to the arguments in \cite{RRS09}.


 When $H$  $\varepsilon$-contains $H_{k}(n,m)$, Lu, Yu and Yuan \cite{LYY} use the structure of $H_{k}(n,m)$ to greedily find the desired matching. (Lemma 2.3 in \cite{LYY})
\begin{lemma}[Lu, Yu and Yuan, \cite{LYY}]\label{LYY-close}
Let $n,m,k$ be integers and $0<\varepsilon<(8^{k-1}k^{5(k-1)}k!)^{-3}$, such that  $k\geq3$, $n\geq 8k^6/(1-5k^2\sqrt{\varepsilon})$, and $n/(2k^4)+1<m\leq n/k$. Let $H$ be a $k$-graph on $n$ vertices such that $\delta_1(H)>{n-1\choose k-1}-{n-m\choose k-1}$ and $H$ $\varepsilon$-contains $H_{k}(n,m)$, then $\nu(H)\geq m$.
\end{lemma}

Thus for completing the proof of Theorem \ref{main-thm}, it suffices to prove the following lemma. By $x\ll y$ we mean that $x$ is sufficiently small compared with $y$ such that $x,y$ satisfy finitely many inequalities in the proof.
\begin{lemma}\label{nonclose}
	Let $k\geq 3 $ be an integer and let $\varepsilon,\rho,\beta$  be constants such that  $0<\varepsilon<(3^{k-2}k!k^3)^{-1}$ and $0<\beta\ll\rho<\varepsilon^4/(18k^2)^4$. Let $n,m$ be two positive integers such that $n$ is sufficiently large, and
	$\frac{n}{k^4}\leq m\leq (\frac{k}{2(k-1)}-\beta)\frac{n}{k}$. Let $H$ be a $k$-graph on $n$ vertices.
	If $H$ does not $\varepsilon$-contain $H_{k}(n,m)$ and $\delta_1(H)>{n-1 \choose k-1}-{n-m\choose k-1}- \rho n^{k}$, then $\nu(H)\geq m$ .
\end{lemma}

In order to prove Lemma \ref{nonclose}, we first construct a $k$-graph $H^k_r$ from $H$ such that $\nu(H)\geq m$ if and only if $H^k_r$ has an almost perfect matching. In Section 2, we follow some ideas from \cite{AFHRRS12,FK18,KOT14,LYY} to prove that $H^k_r$ has a perfect fractional matching by using a stability result proved by Lu, Yu and Yuan \cite{LYY}. In Section 3, we use the two-round randomization method from \cite{AFHRRS12} to show that $H^k_r$ has a nearly regular spanning subhypergraph in which all $2$-degrees are much smaller than the vertex degrees. Then a result of Frankl and R\"odl \cite{FR85} implies that $H^k_r$ has an almost perfect matching.

In section 2, we apply the stability result proved by Lu, Yu and Yuan \cite{LYY} to an $(n-1)$-vertex $(k-1)$-graph. Since their result holds for $m\leq n/(2k)$, we are only able to prove Lemma \ref{ind-fpm} for $m\leq \frac{n-1}{2(k-1)}$. In section 3, we use randomization method and Lemma \ref{ind-fpm} to prove Lemma \ref{nonclose}. In order to prove a random subgraph has a perfect fractional matching by using Lemma \ref{ind-fpm}, we need $m\leq (\frac{k}{2(k-1)}-\beta)\frac{n}{k}$, where $\beta$ is a small positive.

We end this section with additional notations. For any positive integer $n$, let $[n]:=\{1,\ldots,n\}$. For a  $k$-graph $H$ and  $S\subseteq V(H)$, we use $H-S$ to denote the hypergraph obtained from $H$ by deleting $S$ and all edges of $H$ intersecting set $S$, and we use $H[S]$ to denote the subhypergraph with vertex set $S$ and edge set $\{e\in E(H) : e\subseteq S\}$. For a $k$-graph $H$ and a vertex $v\in V(H)$, let $N_H(v):=\{e\in \binom{V(H)}{k-1}: e\cup\{v\}\in E(H)\}$. We omit the floor and ceiling functions when they do not affect the proof.



\section{Fractional Matching}


A {\it fractional matching} in a $k$-graph $H$ is a function $\varphi: E(H)\rightarrow [0,1]$ such that for any $v\in V(H)$, $\sum_{\{e\in E: v\in e\}}\varphi(e)\le 1$.
A fractional matching is called {\it perfect} if $\sum_{e\in E}\varphi (e)=|V(H)|/k$.
For a hypergraph $H$, let $$\nu'(H)=\max\left\{\sum_{e\in E(H)} \varphi(e): \varphi \mbox{ is a fractional matching in $H$}\right\}.$$
A {\it fractional vertex cover} of $H$  is a
function $w:V(H)\rightarrow [0,1]$ such that for each $e\in E$, $\sum_{v\in e}w(v)\ge 1$. Let
$$\tau'(H)=\min\left\{\sum_{v\in V(H)}w(v): w \mbox{ is a fractional vertex cover of $H$}\right\}.$$
Then the strong duality theorem of linear programming gives
$$\nu'(H)=\tau'(H).$$

For a complete $k$-graph, we have the following observation.
\begin{obser}\label{cliquefrac}
Let $n,k$ be two integers such that $n> k\geq 2$. Let $H$ be a $k$-graph on $n$ vertex with edge set $\binom{V(H)}{k}$, then $H$ has a perfect fractional matching.
\end{obser}
\pf
	Let $V(H)=[n]$ and let $e_i=\{i,\ldots,i+k-1\}$ for $i=1,\ldots,n$, where the addition is on modular $n$ (except we write $n$ instead of $0$). Write $E=\{e_1,\ldots,e_{n}\}$. Note that $e_i\in E(H)$. Let $\varphi:E(H)\rightarrow [0,1]$ such that
	\begin{equation*}
		\varphi(e)= \left\{
		\begin{array}{ll}
			1/k, & \hbox{$e\in E$;} \\
			0, & \hbox{otherwise.}
		\end{array}
		\right.
	\end{equation*}
	It is not difficult to see that $\varphi$ is a perfect fractional matching in $H$.
\qed

Recall that $K_r^k$ is a complete $k$-graph with $r$ vertices and edge set ${V(K_r^k)\choose k}$.
Let $H_r^k=H+K_r^k$ denote a $k$-graph with vertex set $V(H)\cup V(K_r^k)$ and edge set
\[
E(H_r^k)=E(H)\cup \left\{ e\in {V(H)\cup V(K_r^k)\choose k} : e\cap V(K_r^k)\neq \emptyset \right\}.
\]

In this section, we prove that $H^k_r$ satisfying the conditions in Lemma \ref{ind-fpm} has a perfect fractional matching. In order to prove $H^k_r$ has a perfect fractional matching, we need the following lemma for stable $k$-graphs (Lemma 4.2 in \cite{LYY}). Let $H$ be a $k$-graph with vertex set $[n]$. For any $\{u_{1},\ldots,u_{k}\},\{v_{1},\ldots,v_{k}\}\in {[n]\choose k}$ with $u_i<u_{i+1}$ and $v_i<v_{i+1}$ for $1\leq i\leq k-1$, we write $\{u_{1},\ldots,u_{k}\}\leq \{v_{1},\ldots,v_{k}\}$ if $u_i\leq v_i$ for $1\leq i\leq k$. $H$ is called \emph{stable} if for $e,f\in \binom{[n]}{k}$ with $e\leq f$, $f\in E(H)$ implies that $e\in E(H)$.
\begin{lemma}[Lu, Yu and Yuan, \cite{LYY}]\label{stafrankl}
Let $k$ be a positive integer, and let $c$ and  $\rho$  be constants such that  $0< c<1/(2k)$ and $0<\rho\leq(1+18(k-1)!/c)^{-2}$.
		Let $n,m$ be positive integers   such that $n$ is sufficiently large and
	$cn\leq m\leq n/(2k)$. Let $H$ be a stable $k$-graph with vertex set $[n]$.
	If $e(H)>{n\choose k}-{n-m\choose k}- \rho n^{k}$ and  $\nu(H)\leq m$, then $H$ $\sqrt{\rho}$-contains $H_{k}([n]\setminus[m],[m])$\footnote{In \cite{LYY}, the conclusion is $H$ $\sqrt{\rho}$-contains $H_{k}^k([n]\setminus[m],[m])$, where $H_{k}^k([n]\setminus[m],[m])$ is the $k$-graph with vertex set $[n]$ and edge set $\{e\in \binom{[n]}{k}:e\cap [m]\neq \emptyset\}$. Since $H_{k}([n]\setminus[m],[m])$ is a subgraph of $H_{k}^k([n]\setminus[m],[m])$, this conclusion implies $H$ $\sqrt{\rho}$-contains $H_{k}([n]\setminus[m],[m])$}.
\end{lemma}



Let $H$ be a $k$-graph on $n$ vertices and $U,W$ be a partition of $V(H)$. Given $0<\theta<1$, a vertex $v\in V(H)$ is \emph{$\theta$-good} with respect to $H_{k}(U,W)$ if $|N_{H_{k}(U,W)} (v)\setminus N_H(v)|\leq \theta n^{k-1}$. Otherwise 
$v$ is \emph{$\theta$-bad}. A set $I \subseteq V(H)$ that contains no edge of $H$ is called an
\emph{independent set}  in $H$. We use $\alpha(H)$ to denote the size of the largest independent set in the hypergraph $H$.

\begin{lemma}\label{ind-fpm}
Let $k\geq 3$ be an integer and let $\rho,\varepsilon$ be constants such that $0<\varepsilon\leq (3^{k-2}k!k^3)^{-1}$ and $0<\rho<\varepsilon^4/(2k^8)$. Let $n,m,r$ be integers such that $n$ is sufficiently large, $\frac{n-1}{2k^4}\leq m\leq \frac{n-1}{2(k-1)}$ and $(r-k)(k-1)\geq n-km$. Let $H$ be a $k$-graph on $n$ vertices such that $\alpha(H)<n-m-\varepsilon n$. If $\delta_1(H)>{n-1\choose k-1}-{n-m\choose k-1}- \rho n^{k-1}$, then $H_r^k$ has a  perfect fractional  matching.
\end{lemma}

\pf 
Let $V(H)=[n]$ and let $Q=V(H_r^k)\setminus V(H)=\{n+1,\ldots,n+r\}$.  Let $w:V(H_r^k)\rightarrow [0,1]$ be a minimum fractional vertex cover of  $H_r^k$. Rename the vertices in $[n]$ such that
\begin{equation}\label{weight}
w(1)\geq w(2)\geq\cdots\geq w(n).
\end{equation}
Let $H'$ be a $k$-graph with vertex set $V(H_r^k)$ and edge set

\begin{align*}
 E(H')=\left\{e\in {V(H_r^k)\choose k} : \sum_{x\in e}w(x)\geq 1\right\}.
\end{align*}
Since $\sum_{x\in e}w(x)\geq 1$ for every $e\in E(H^k_r)$, $H'$ is a superhypergraph of $H^k_r$. Let $G=H'-Q$. Thus $G$ is a superhypergraph of $H$. Let $G'$ be the $(k-1)$-graph with vertex set $[n-1]$ and edge set $N_G(n)$.

\begin{cla}
	
\begin{enumerate}[itemsep=0pt,parsep=0pt,label=$($\roman*$)$]
	\item $G'$ is stable.
	\item Let $S:=[m+\varepsilon n]$, then $G[S]$ is a complete $k$-graph.
	\item For any $e\in N_G(n)$, if $i\in[n]\setminus e$, then $e\in N_G(i)$.
\end{enumerate}
\end{cla}
\pf
For two sets $\{x_1,\ldots,x_{k-1}\},\{y_1,\ldots,y_{k-1}\}\in \binom{[n-1]}{k-1}$ such that
$x_i\leq y_i$ for $1\leq i\leq k-1$, if $\{y_1,\ldots,y_{k-1}\}\in E(G')$, then $\{y_1,\ldots,y_{k-1},n\}\in E(G)$. Thus $\sum_{i=1}^{k-1}w(y_i)+w(n)\geq 1$. By (\ref{weight}), we can derive that $\sum_{i=1}^kw(x_i)+w(n)\geq 1$. That is, $\{x_1,\ldots,x_{k-1}\}\in E(G')$. Thus $G'$ is stable, and (i) follows.  To prove (ii), suppose that $G[S]$ is not a complete $k$-graph. Then there exists a set $\{v_1,\ldots,v_k\}\in \binom{S}{k}$ such that $\{v_1,\ldots,v_k\}\notin E(H')$. Thus $\sum_{i=1}^kw(v_i)< 1$. Let $S':=\{m+\varepsilon n-1,\ldots,n\}$. By (\ref{weight}), for every $e\in \binom{S'}{k}$, we have $\sum_{x\in e}w(x)\leq \sum_{i=1}^kw(v_i)< 1$. Thus $S'$ is an independent set in $G$ with $|S'|\geq n-m-\varepsilon n$, contradicting the fact that $ \alpha(G)\leq \alpha(H)<n-m-\varepsilon n$. To prove (iii), let $e\in N_G(n)$. Since $E(G)\subseteq E(H')$, we have $\sum_{x\in e}w(x)+w(n)\geq 1$. Thus by (\ref{weight}), $\sum_{x\in e}w(x)+w(i)\geq 1$ for any $i\in[n]\setminus e$. That is, $e\in N_G(i)$ for any $i\in[n]\setminus e$.
\qed

One can see that $w$ is also a fractional vertex cover of $H'$. Thus $w$ is also a minimum fractional vertex cover of $H'$. By Linear Programming Duality Theory, we have
$\nu'(H_r^k)=\tau'(H_r^k)=\tau'(H')=\nu'(H')$. Thus it suffices to prove that $H'$ has a perfect fractional matching.
Suppose that $n+r\equiv s\pmod k$, where $0\leq s\leq k-1$.  Let $Q'=\{n+1,\ldots,n+s\}$. We first find a perfect matching in $H'-Q'$, then use it to construct a  perfect fractional matching in $H'$.

\begin{cla}
	$\nu(G)\geq m$.
\end{cla}


\pf  Let $W:=[m]$, $U:=[n-1]\setminus[m]$. Assume $G'$ does not $\sqrt{2\rho}$-contain $H_{k-1}(U,W)$. Let $c=1/(2k^4)$ be a constant as in Lemma \ref{stafrankl}. Note that $2\rho<\varepsilon^4/k^8< (k!k^5)^{-4}\leq (1+18(k-2)!/c)^{-2}$, $\frac{n-1}{2k^4}\leq m\leq\frac{n-1}{2(k-1)}$ and $e(G')=|N_G(n)|> {n-1\choose k-1}-{n-1-m\choose k-1}- 2\rho (n-1)^{k-1}$. By Claim 1, $G'$ is stable. We can derive that $\nu(G')> m$ by the contrapositive of Lemma \ref{stafrankl}. Recall that $n\geq2(k-1)m+1$ and $k\geq 3$. Thus we have $n\geq km$. Let $M_1=\{f_1,\ldots,f_{m}\}$ be a matching of size $m$ in $G'$ and let $\{v_{1},\ldots,v_{{m}}\}\subseteq [n]\setminus\left( \bigcup_{i=1}^{m} f_i\right)$.
By  Claim 1(iii) ,  we have $M_1\subseteq N_G(v)$ for any $v\in [n]\setminus\left( \bigcup_{i=1}^{m} f_i\right)$. Thus $M_{1}'=\{f_i\cup \{v_{i}\} : i\in [m]\}$ is a matching of size $m$ in $G$.

So we may assume that $G'$ $\sqrt{2\rho}$-contains $H_{k-1}(U,W)$. Then $G'$ contains less than $(k-1)(2\rho)^{1/4}n$ $(2\rho)^{1/4}$-bad vertices with respect to $H_{k-1}(U,W)$. Otherwise,
 \begin{align*}
 |E(H_{k-1}(U,W)\setminus E(G')|\geq &\frac{1}{k-1}\sum_{v\in V(G')}|N_{H_{k-1}(U,W)}(v)\setminus N_{G'}(v)|\\
 >&\frac{1}{k-1}\cdot (k-1)(2\rho)^{1/4}n\cdot(2\rho)^{1/4}(n-1)^{k-2}\\
 >&\sqrt{2\rho} (n-1)^{k-1},
 \end{align*}
 a contradiction.
 Let $B$ denote the set of $(2\rho)^{1/4}$-bad vertices in $W$. Write $b:=|B|$.
 So  $b<(k-1)(2\rho)^{1/4}n$.

First we find a matching $M_{21}$ of size $b$ in $G[U]$. Let $S_1:=\{m,\ldots,m+\varepsilon n\}$. By  Claim 1(ii), we can derive that $G[S_1]$ is a complete $k$-graph.  Since $b< (k-1)(2\rho)^{1/4}n$ and $2\rho<\varepsilon^4/k^8$, we can find pairwise disjoint edges $f_1,\ldots,f_{b}$ in $G[S_1]$. Thus $M_{21}=\{f_1,\ldots,f_{b}\}$ is a matching in $G[U]$.

Let $U_1:=U\setminus V(M_{21})$, $W_1:=W\setminus B$ and $G'':=G'-(V(M_{21})\cup B)$. Thus $V(G'')=U_1\cup W_1$ and $|W_1|=m-b$. For every $x\in W_1$, since $x$ is $(2\rho)^{1/4}$-good in $G'$ with respect to $H_{k-1}(U,W)$ and $N_{H_{k-1}(U_1,W_1)}(x)\setminus N_{G''}(x)\subseteq N_{H_{k-1}(U,W)}(x)\setminus N_{G'}(x)$, we have
\begin{align*}
	\left|N_{H_{k-1}(U_1,W_1)}(x)\setminus N_{G''}(x)\right|
	\leq \left|N_{H_{k-1}(U,W)}(x)\setminus N_{G'}(x)\right|\leq (2\rho)^{1/4}(n-1)^{k-2}<\binom{n/3}{k-2}. 
\end{align*}
It follows that
\begin{align}
\left|N_{G''}(x)\cap \binom{U_1}{k-2}\right|\geq &\left|N_{H_{k-1}(U_1,W_1)}(x)\cap \binom{U_1}{k-2}\right|-\left|N_{H_{k-1}(U_1,W_1)}(x)\setminus N_{G''}(x)\right|\notag\\
>& \binom{|U_1|}{k-2}-\binom{n/3}{k-2}.\label{large-degree}
\end{align}
That is, every vertex $x\in W_1$ has large degree in $G''$.

Now we use vertices in $W_1$ to construct a matching $M_{22}=\{e_1,\ldots,e_{m-b}\}$ in $G''$ such that $|e_i\cap W_1|=1$ for $i=1,\dots,m-b$. Suppose for some integer $0< t\leq m-b-1$, we have found a matching $\{e_1,\ldots,e_t\}$ in $G''$ such that $|e_i\cap W_1|=1$ for $1\leq i\leq t$. Note that
\begin{align}
\left|U_1\setminus\left(\bigcup_{i=1}^t e_i\right)\right|&=|U|-|V(M_{21})|-\left|\left(\bigcup_{i=1}^t e_i\right)\cap U_1\right|\notag\\
&>n-1-m-kb-(k-2)(m-b)\notag\\
&=n-(k-1)m-2b-1\notag\\
&\geq n-(k-1)m-2(k-1)(2\rho)^{1/4}n-1\quad\mbox{(since $b< (k-1)(2\rho)^{1/4}n$)}\notag\\
&>n/3 \quad \mbox{(since $m\leq (n-1)/(2k-2)$ and $2\rho<\varepsilon^4/k^8<1/k^{20}$)}.\label{manyvertices}
\end{align}
Let $x\in W_1\setminus \left(\bigcup_{i=1}^t e_i\right)$, by inequalities (\ref{large-degree}) and (\ref{manyvertices}), we have
\begin{align*}
\left|N_{G''}(x)\cap\binom{U_1\setminus(\bigcup_{i=1}^t e_i)}{k-2}\right|\geq \left|N_{G''}(x)\cap \binom{U_1}{k-2}\right|-\left(\binom{|U_1|}{k-2}-\binom{|U_1\setminus(\bigcup_{i=1}^t e_i)|}{k-2}\right)>0.
\end{align*}
Thus there exists an edge $e_{t+1}\subseteq V(G'')\setminus (\bigcup_{i=1}^t e_i)$ such that $e_{t+1}\cap W_1=\{x\}$. Continue this process until $t=m-b-1$. Then $M_{22}=\{e_1,\ldots,e_{m-b}\}$ is the desired matching.

Recall that $n\geq km$. Let $v_{1},\ldots, v_{{m-b}}$ be $m-b$ vertices of $G-V(M_{21}\cup M_{22})$. By  Claim 1(iii), we have $M_{22}\subseteq N_{G}(v_{i})$ for $1\leq i\leq m-b$. Thus
$M_{22}'=\{e_i\cup \{v_{i}\} : i\in [m-b]\}$ is a matching of size $m-b$ in $G-V(M_{21})$. Then $M_{21}\cup M_{22}'$ is a matching of size $m$ in $G$ and thus $\nu(G)\geq m$.
\qed

Let $M$ be a matching of size $m$ of $G$.
Note that $N_{H'}(n+i)={[n+r]\setminus \{n+i\}\choose k-1}$ for $1\leq i\leq r$ and
\[
r-s\geq (n-km)/(k-1).
\]
Thus $H'-Q'-V(M)$ has a perfect matching, say $M'$. For the case $s=0$, $M\cup M'$ is a perfect matching in $H'$. For the case $s\neq 0$, let $f\in M'$. By Observation \ref{cliquefrac}, $H'[f\cup Q']$ has a  perfect fractional matching $\varphi$. Let $\varphi' :E(H')\rightarrow [0,1]$ such that
\begin{equation*}
  \varphi'(f)=\left\{
     \begin{array}{ll}
       1, & \hbox{$e\in M\cup (M'-f)$;} \\
       \varphi(e), & \hbox{$e\in E(H'[f\cup Q'])$;} \\
       0, & \hbox{otherwise.}
     \end{array}
   \right.
\end{equation*}
Recall that $M\cup (M'-f)$ is a perfect matching in $H'-(Q'\cup f)$. So $\varphi'$ is a  perfect fractional matching in $H'$. This completes the proof. \qed

%
%

%

\section{Almost Perfect Matching}



The following lemma asserts that the existence of an almost perfect matching in any nearly regular $k$-graph in which all $2$-degrees are much smaller than the vertex degrees (see Theorem 1.1 in \cite{FR85} or Lemma 4.2 in \cite{AFHRRS12}). For any positive integer $\ell$, we use $\Delta_\ell(H)$ to denote the maximum $\ell$-degree of a hypergraph $H$.
\begin{lemma}[Frankl and R\"odl, \cite{FR85}]\label{FR-1d2d}
	For every integer $k \geq 2$ and any real $\sigma>0$,
	there exist $\tau=\tau(k,\sigma)$,
	$d_0=d_0(k,\sigma)$ such that
	for every $n \ge D \ge d_0$ the following holds:
	Every $k$-graph $H$ on $n$ vertices with $(1-\tau)D <d_H(v)<(1+\tau)D$ and $\Delta_2(H) <\tau D$ contains a matching covering all but at most $\sigma n$ vertices.
\end{lemma}

Let $Bi(n, p)$ be the binomial distribution with parameters $n$ and $p$.
The following lemma on Chernoff bound can be found in Alon and
Spencer \cite{AS08} (page 313).

\begin{lemma}[Chernoff]\label{chernoff}
Suppose $X_1,\ldots, X_n$ are independent random variables taking values in $\{0, 1\}$. Let $X=\sum_{i=1}^n X_i$ and $\mu = \mathbb{E}[X]$. Then, for any $0 < \delta \leq 1$,
$$\mathbb{P}[X \geq (1+\delta) \mu] \le  e^{-\delta^2\mu/3}\mbox{ and } \mathbb{P}[X \leq (1-\delta) \mu] \le  e^{-\delta^2 \mu/2}.$$
In particular, when $X \sim Bi(n,p)$ and $\lambda<\frac{3}{2}np$,  then
\[
\mathbb{P}(|X-np|\geq \lambda)\leq e^{-\Omega(\lambda^2/np)}.\label{cher}
\]
\end{lemma}

In order to find a spanning subgraph in a hypergraph satisfying conditions in Lemma~\ref{FR-1d2d}, we use the same two-round randomization
technique as in \cite{AFHRRS12}. The only difference is that between
the two rounds, we also need to bound the independence number of the subgraph. 
%
The following lemma (Lemma 5.4 in \cite{LYY}) was proved by Lu, Yu and Yuan using hypergraph container method.
\begin{lemma}[Lu, Yu and Yuan, \cite{LYY}]\label{indep}
Let $c, \theta, \zeta$ be positive reals and let $k, n$ be positive integers. Let $H$ be an $n$-vertex $k$-graph such that $e(H)\ge cn^k$  and $e(H[S])\ge\theta e(H)$ for all $S\subseteq V(H)$ with $|S|\ge \zeta n$. 
Let $R\subseteq V(H)$ be obtained  by taking each vertex of $H$ uniformly at random with probability $n^{-0.9}$.
Then for any real number $\gamma$ with $0<\gamma\ll \zeta $, there exists an integer $n_0$ such that  for any $n>n_0$, $\alpha(H[R])\le (\zeta +\gamma)n^{0.1}$ 
with probability at least $1-e^{-\Omega (n^{0.1})}$.
\end{lemma}

\begin{lemma}\label{Close-Inde}
	Let $n,k,m$ be integers such that $k\geq 3$,  $n$ is sufficiently large and $n/2k^4\leq m< n/k$. Let $0<\epsilon<1/k$ and $0<\varrho< \epsilon/12$.  Let $H$ be a $k$-graph on $n$ vertices. If $\delta_1(H)\geq {n-1\choose k-1}-{n-m\choose k-1}-\varrho n^{k-1}$ and $H$ does not $\epsilon$-contain $H_k(n,m)$, then $e(H[S])\geq \frac{\epsilon n^k}{2k^2}$ for any set $S\subseteq V(H)$ with $|S|\geq (1-\frac{m}{n}-\frac{\epsilon}{7})n$.
\end{lemma}

\pf Suppose that the result does not hold. Then $H$ has a set $A$ such that $|A|\geq (1-\frac{m}{n}-\frac{\epsilon}{7})n$ and $e(H[A])<\frac{\epsilon n^k}{2k^2}$. After removing vertices from $A$ if necessary, we may assume that $|A|\leq n-m$. Let $ W\subseteq V(H)\setminus A$ such that $|W|=m-1$. Let $U=V(H)\setminus W$ and $B=U\setminus A$.
One can see that $|B|\leq \epsilon n/7+1$.
Since $e(H[A])<\frac{\epsilon n^k}{2k^2}$, the number of edges belonging to $H[U]$ is at most
\begin{align*}
\sum^k_{i=1}\binom{\left|B\right|}{i}\binom{\left|U\setminus B\right|}{k-i}+e(H[A])&<\sum^k_{i=1}\frac{(\frac{\epsilon n}{7}+1)^in^{k-i}}{i!(k-i)!}+\frac{\epsilon n^k}{2k^2}\\
&\leq\sum^k_{i=1}\frac{\epsilon^i n^k}{6^ii!(k-i)!}+\frac{\epsilon n^k}{2k^2}\\
&\leq\frac{k\epsilon n^k}{6(k-1)!}+\frac{\epsilon n^k}{2k^2}
\end{align*}
%
%
for sufficiently large $n$. So we have
\begin{equation*}
\sum_{x\in U}\left|N_{H[U]}(x)\right|= k\cdot e(H[U])\leq k\left(\frac{k\epsilon n^k}{6(k-1)!}+\frac{\epsilon n^k}{2k^2}\right).
\end{equation*}
Since every edge of $H_k(U,W)$ intersects $U$, we may infer that
\begin{align*}
	&\left|E(H_{k}(U,W)) \setminus E(H)\right|\leq \sum_{x\in U}\left|N_{H_{k}(U,W)}(x)\setminus N_{H}(x)\right|\\
	= &\sum_{x\in U}(\left|N_{H_{k}(U,W)}(x)\right|-\left|N_{H}(x)\setminus N_{H[U]}(x)\right|)\\
	\leq &\sum_{x\in U}\left({n-1\choose k-1}-{n-m\choose k-1}-\left({n-1\choose k-1}-{n-m\choose k-1}-\varrho n^{k-1}-\left|N_{H[U]}(x)\right|\right)\right)\\
	= &\left|U\right|\left({n-1\choose k-1}-{n-m\choose k-1}-\left({n-1\choose k-1}-{n-m\choose k-1}-\varrho n^{k-1}\right)\right)+\sum_{x\in U}\left|N_{H[U]}(x)\right|\\
	\leq &\left|U\right|\varrho n^{k-1}+k(\frac{k\epsilon n^k}{6(k-1)!}+\frac{\epsilon n^k}{2k^2})\\
	\leq&\left|U\right|\varrho n^{k-1}+3\epsilon n^k/4+\epsilon n^k/2k\quad \mbox{(since $k\geq 3$)}\\
	\leq& \epsilon n^k,
\end{align*}
where the second equality follows from $N_{H_{k}(U,W)}(x)\cap N_H(x)=N_{H}(x)\setminus N_{H[U]}(x)$ for each $x\in U$.
Hence $H$ $\epsilon$-contains $H_{k}(U,W)$, a contradiction. \qed

The following lemma can be found in \cite{AFHRRS12,LYY} (Lemma 5.5 in \cite{LYY}), which is the first round of randomization.
\begin{lemma}\label{Random}
	Let $k\geq 3$ be an integer. Let $H$ be a $k$-graph on $n$ vertices. Take $n^{1.1}$ independent copies of $R$ and denote them by $R^i$, $1\le i\le n^{1.1}$, where $R$ is chosen from $V(H)$ by taking each vertex uniformly at random with probability $n^{-0.9}$ and then deleting less than $k$ vertices uniformly at random so that $|R|\in k\mathbb{Z}$.
    For each $S\subseteq V(H)$, let $Y_S:=|\{i: \ S\subseteq R^i\}|$.
    Then with probability at least $1-o(1)$,  all of the following statements hold:
        \begin{itemize}
              \setlength{\itemsep}{0pt}
		\setlength{\parsep}{0pt}
		\setlength{\parskip}{0pt}
            \item [$($i$)$] for every $v\in V$, $Y_{\{v\}}=(1+o(1)) n^{0.2}$,
            \item [$($ii$)$] every pair $\{u, v\} \subseteq V$ is contained in at most two copies $R^i$,
            \item [$($iii$)$] every edge $e \in E(H)$ is contained in at most one copy $R^i$,
            \item [$($iv$)$] for all $i= 1, \dots ,n ^{1.1}$, we have $||R^i|-n^{0.1}|\leq n^{0.095}$, and
            \item [$($v$)$]  if $\mu,\rho$ are constants with $0<\mu\ll\rho$, $n/k-\mu n\leq m\leq n/k$, and $\delta_1(H)\geq\binom{n-1}{k-1}-\binom{n-m}{k-1}-\rho n^{k-1}$, then for all $i=1,\ldots,n^{1.1}$ and any positive real $\rho'\geq 2\rho$, we have
            \[
            \delta_1(H[R_i])>\binom{|R^i|-1}{k-1}-\binom{|R^i|-|R^i|/k}{k-1}-\rho'|R^i|^{k-1}.
            \]


        \end{itemize}
    \end{lemma}

We summarize the second round of randomization in \cite{AFHRRS12} as the following lemma (see the proof of Claim 4.1 in \cite{AFHRRS12}).
\begin{lemma}\label{secondrandom}
Assume $R^i$, $i=1,2,\ldots,n^{1.1}$ satisfy $(i)-(v)$ in Lemma \ref{Random}, and each $H[R^i]$ has a perfect fractional  matching $\varphi^i$. Then there exists a spanning subgraph $H'$ of $H$ such that $d_{H'}(v)=(1+o(1)) n^{0.2}$ for each $v\in V$, and $\Delta_2(H')\leq n^{0.1}$.
\end{lemma}

\noindent \textbf{Proof of Lemma \ref{nonclose}.} 
Let $V(H)=[n]$ and $\eta=\beta/3$.
We choose an integer $r$ such that
\begin{align}\label{GL1}
r=\left\lceil\frac{n-km-\eta n}{k-1}\right\rceil.
\end{align} Let $Q:=V(K_r^k)=\{n+1,\ldots,n+r\}$. Recall that $H_r^k=H+K_r^k$ and let $n_1:=n+r$.
By $\delta_1(H)>\binom{n-1}{k-1}-\binom{n-m}{k-1}-\rho n^{k-1}$, we can derive that
\begin{equation}\label{GL2}
\delta_1(H_r^k)>{n_1-1\choose k-1}-{n_1-m-r\choose k-1}-\rho n^{k-1}.
\end{equation}
It suffices to show that $\nu(H_r^k)\geq m+r$. Indeed, if there is a matching $M$ of size $m+r$ in $H_r^k$, then there are at most $r$ edges in $M$ intersecting $Q$ and thus $\nu(H)\geq m$.

Let $R\subseteq V(H_r^k)$ be obtained  by taking each vertex of $H_r^k$ uniformly at random with probability $n_1^{-0.9}$.  Take $n_1^{1.1}$ independent copies of $R$ and denote them by $R^i$, $1\le i\le n_1^{1.1}$.

By Lemma~\ref{Random}(iv), we have
\begin{align}\label{GL3}
	n_1^{0.1}-n_1^{0.095}\leq |R^i|\leq n_1^{0.1}+n_1^{0.095}\mbox{ for all $i=1,\ldots, n_1^{1.1}$}
\end{align}
with probability $1-o(1)$. One can see that $|V(H)|\geq\frac{k-1}{k}|V(H^k_r)|\geq \frac{2}{3}|V(H^k_r)|$ since $r(k-1)<n$ and $k\geq 3$. For each $i$, $|R^i\cap V(H)|$ is a binomial random variable with expectation $nn_1^{-0.9}$. Applying Lemma \ref{cher} with $\lambda=n_1^{0.095}$, we have
\begin{align*}
\mathbb{P}(\left|\left|R^i\cap V(H)\right|-nn_1^{-0.9}\right|\geq n_1^{0.095})\leq e^{-\Omega(n_1^{0.09})}.
\end{align*}
Thus by the union bound, we have
\begin{align}\label{GL4}
nn_1^{-0.9}-n_1^{0.095}\leq \left|R^i\cap V(H)\right|\leq nn_1^{-0.9}+n_1^{0.095} \mbox{ for all $i=1,\ldots, n_1^{1.1}$}
\end{align}
 with probability at least $1-n_1^{1.1}e^{-\Omega(n_1^{0.09})}$. Write $r_i:=|R^i\cap Q|$.
With similar discussion, one can see that
\begin{align}\label{GL5}
rn_1^{-0.9}-n_1^{0.095}\leq r_i\leq rn_1^{-0.9}+n_1^{0.095} \mbox{ for all $i=1,\ldots, n_1^{1.1}$}
\end{align} with probability at least $1-n_1^{1.1}e^{-\Omega(n_1^{0.09})}$.
Thus by (\ref{GL1}), (\ref{GL4}) and (\ref{GL5}), we have
\begin{align}\label{cond:2}
(r_i-k)(k-1)&\geq(rn_1^{-0.9}-n_1^{-0.095}-k)(k-1)\notag\\
&\geq (n-km-\eta n)n_1^{-0.9}-(n_1^{0.095}+k)(k-1)\notag\\
&\geq \left|R^i\cap V(H)\right|-kmn_1^{-0.9}-2\eta nn_1^{-0.9} \mbox{ for all $i=1,\ldots, n_1^{1.1}$}
\end{align}
with probability $1-o(1)$.

Since $H$ does not $\varepsilon$-contain $H_k(n,m)$, then 
 by Lemma~\ref{Close-Inde},
$e(H[S])\ge \varepsilon n^k/2k^2$
for all $S\subseteq V(H)$ with $|S|\ge \alpha n$, where $\alpha=1-m/n-\varepsilon/7$. Since each vertex in $Q$ has degree $\binom{n_1-1}{k-1}$, replacing a vertex in $S$ by a vertex in $Q$ will not decrease the number of edges in $H^k_r[S]$. Thus for every $S\subseteq V(H_r^k)$ with $|S|\ge \alpha n$, we have $e(H^k_r[S])\ge \varepsilon n^k/2k^2\ge (\frac{k-1}{k})^k\varepsilon n_1^k/2k^2$. Then by Lemma~\ref{indep}, with probability $1-o(1)$, for each $i$, $H_r^k[R^i]$ has  independence number $\alpha(H_r^k[R^i])\le (n_1-m-r-\varepsilon n/8)n_1^{-0.9}=(n-m-\varepsilon n/8)n_1^{-0.9}$. Note that $\alpha(H[R^i\cap V(H)])=\alpha(H_r^k[R^i])$. 
  So we have \begin{align}\label{cond:1}
    \alpha(H[R^i\cap V(H)])\le (1-m/n-2\eta/k-\varepsilon/9)\left|R^i\cap V(H)\right| \mbox{ for all $i=1,\ldots, n_1^{1.1}$}
    \end{align} with probability $1-o(1)$.

    Note that  $n_1/k-2\eta n_1\leq m+r\leq n_1/k$, where $2\eta\ll\rho$. By  Lemma~\ref{Random}(v) and inequality (\ref{GL2}), with probability $1-o(1)$, for all $i=1,\ldots, n^{1.1}$, we have
    \begin{align}\label{random-degree}
    \delta_1(H_r^k[R^i])> \binom{|R^i|-1}{k-1}-{|R^i|-|R^i|/k\choose k-1}-3\rho|R^i\cap V(H)|^{k-1}.
    \end{align}
     So by inequalities (\ref{GL3}), (\ref{GL5}), (\ref{random-degree}) and $n_1/k-2\eta n_1\leq m+r\leq n_1/k$, with probability $1-o(1)$, for all $i=1,\ldots, n_1^{1.1}$, we have
     \begin{align}
    &\delta_1(H[R^i\cap V(H)])=\delta_1(H_r^k[R^i])-\left(\binom{\left|R^i\right|-1}{k-1}-\binom{\left|R^i\cap V(H)\right|-1}{k-1}\right)\notag\\
    &> \binom{\left|R^i\cap V(H)\right|-1}{k-1}-{\left|R^i\right|-\left|R^i\right|/k\choose k-1}-3\rho\left|R^i\cap V(H)\right|^{k-1}\notag\\
    &=\binom{\left|R^i\cap V(H)\right|-1}{k-1}-{\left|R^i\cap V(H)\right|-(\left|R^i\right|-kr_i)/k\choose k-1}-3\rho\left|R^i\cap V(H)\right|^{k-1}\notag\\
    &>\binom{|R^i\cap V(H)|-1}{k-1}-{|R^i\cap V(H)|-(m+2\eta n/k)n_1^{-0.9}\choose k-1}-4\rho|R^i\cap V(H)|^{k-1}.\label{cond:3}
    \end{align}
By $\frac{n}{k^4}\leq m\leq (\frac{k}{2(k-1)}-\beta)\frac{n}{k}$ and inequality (\ref{GL4}), with probability $1-o(1)$, we have $\frac{|R^i\cap V(H)|-1}{k^4}\leq(m+2\eta n/k)n_1^{-0.9}\leq \frac{|R^i\cap V(H)|-1}{2(k-1)}$  for all $i=1,\ldots, n_1^{1.1}$.
Hence by Lemma~\ref{ind-fpm} and   by (\ref{cond:2}), (\ref{cond:1}) and $(\ref{cond:3})$, with probability $1-o(1)$, for all  $i=1,\ldots, n_1^{1.1}$, $H_r^k[R^i]$ has a perfect fractional  matching.

Now for the $k$-graph $H^k_r$, we have chosen  $n_1^{1.1}$ subgraphs  $R^1,\ldots, R^{n_1^{1.1}}$ such that (i)-(v) in Lemma \ref{Random} hold and   $H_r^k[R^i]$ has a perfect fractional  matching for  $1\leq i\leq n_1^{1.1}$. Then by Lemma \ref{secondrandom},  there is a spanning subgraph $H'$ of $H_r^k$ such that $d_{H'}(v)=(1+o(1)) n_1^{0.2}$ for any vertex $v$, and $\Delta_2(H')\leq n_1^{0.1}$.

Thus we may apply Lemma~\ref{FR-1d2d} to find a matching covering all but at most $\sigma n_1$ vertices in $H_r^k$, where $\sigma<2\eta/3$ is a positive constant. Hence we have $\nu(H_r^k)\geq (n+r-\sigma n_1)/k>m+r$ by (\ref{GL1}).
This completes the proof.
\qed

\section{Proof of Theorem \ref{main-thm}}
Since the case $n\geq 3k^2m$ was proved by Huang and Zhao \cite{HZ17}, we may assume $m\geq n/3k^2$ in the proof.
let $\rho,\varepsilon$ be constants such that $0<\varepsilon<(8^{k-1}k^{5(k-1)}k!)^{-3}< (3^{k-2}k!k^3)^{-1}$ and $0<\rho<\varepsilon^4/(18k^2)^4$.
For the case $H$ does not $\varepsilon$-contain $H_{k}(n,m)$, let $\beta_0\leq\beta$ be a constant such that $0<\beta_0\ll \rho$. By Lemma \ref{nonclose}, $\nu(H)\geq m$ for $n/3k^2\leq m\leq (\frac{k}{2(k-1)}-\beta_0)n/k$ and sufficiently large $n$. For the case $H$ $\varepsilon$-contains $H_{k}(n,m)$, $\nu(H)\geq m$ for $n/3k^2\leq m\leq n/k$ and sufficiently large $n$ by Theorem \ref{LYY-close}.\qed

\medskip \noindent
{\bf Acknowledgment.} The authors would like to thank Dr. Hao Huang for his valuable suggestions and comments.

\end{document}